\documentclass[10pt,a4paper]{article}

\usepackage{graphicx}
\usepackage{amsmath,amsthm}
\usepackage[numbers]{natbib}
\usepackage{microtype}
\usepackage[hmargin=3cm,vmargin={3.5cm,4cm}]{geometry}
\usepackage[bitstream-charter]{mathdesign}
\usepackage[final]{showkeys}

\numberwithin{equation}{section}
\allowdisplaybreaks

\newtheorem{Theorem}{Theorem}[section]
\newtheorem{Remark}{Remark}[section]
\newtheorem{Definition}{Definition}[section]

\title{\textbf{Coupled systems of fractional equations related to sound propagation: analysis and discussion}}

\author{$\text{Roberto Garra}_1$, $\text{Federico Polito}_2$\\
	\footnotesize (1) -- Dipartimento di Scienze di Base e Applicate per l'Ingegneria,\\
	\footnotesize ``Sapienza'' University, Rome.\\
	\footnotesize Via A. Scarpa 16, 00161, Rome, Italy.\\
	\footnotesize Email address: rolinipame@yahoo.it \\
	\footnotesize (2) -- Dipartimento di Matematica,\\
	\footnotesize University of Torino.\\
	\footnotesize Via Carlo Alberto 10, 10123, Torino, Italy.\\
	\footnotesize Tel: +39-011-6702937, fax: +39-011-6702878\\
	\footnotesize Email address: federico.polito@unito.it
	}

\begin{document}

	\maketitle
	
	\begin{abstract}
		In this note we analyse the propagation of a small density perturbation in a one-dimensional compressible
		fluid by means of fractional calculus modelling, replacing thus the ordinary time derivative with the Caputo
		fractional derivative in the constitutive equations.
		By doing so,
		we embrace a vast phenomenology, including subdiffusive, superdiffusive and also
		memoryless processes like classical diffusions.
		From a mathematical point of view, we study systems of coupled fractional equations,
		leading	to fractional diffusion equations or to equations with sequential fractional derivatives.
		In this framework we also propose a method to solve partial differential equations
		with sequential fractional derivatives by analysing the corresponding
		coupled system of equations.
		
		\smallskip
		
	\end{abstract}

	\section{Introduction}
 	
		Fractional calculus modelling has recently had different developments in the applied science \cite{book}.
		The Caputo derivative of real order introduces 
		a memory formalism as it is an integro-differential operator defined by the convolution of the ordinary
		derivative with a power law kernel.   
		In this note we propose an analysis of the propagation of a small density perturbation through a compressible
		fluid using techniques related to fractional calculus. In practice we replace
		the ordinary time derivative with the Caputo fractional derivative in the constitutive equations
		aiming at furnishing the model with a memory.
		This approach leads to a coupled system of fractional differential
		equations that has to be properly treated.
		
		From a physical point of view, this means introducing a memory formalism in a simple model
		of sound propagation through a one-dimensional compressible fluid. 
		By means of this model we can take into account a time delay 
		in the wave propagation. Similar fractional generalisations in modelling sound
		propagation were already considered in some papers.
		For instance in some recent works, Fellah et al.\ \cite{Fellah1,Fellah2} studied transient-wave
		propagation in porous materials using fractional modelling to take into account
		the frequency variability of some dynamic cofficients of the medium
		like tortuosity and compressibility. Their approach is also validated by experimental observations
		\cite{Fellah3}.
		
		Furthermore, we have also to notice that Tarasov \cite{Tarasov} gave a space-fractional
		formulation of the hydrodynamic equations to describe fluid flow in fractal media. He also studied sound
		wave propagation within this conceptual framework. 
		In our formulation we use a time-fractional generalisation of the characteristic master equations;
		therefore our model is in some way similar to the one of Fellah et al.\ \cite{Fellah2} 
		concerning propagation in porous media, but more general, as we are going to discuss.  
		As a matter of fact we want to show that by using our approach, with some specific hypotheses
		on the range of variability of the fractional order of differentiation,
		a fractional diffusion-wave equation will result by decoupling the constitutive equations.
		This result appears interesting because we obtain a family of equations, depending on
		the real order of differentiation, which includes the heat equation as a particular case.
		We find thus that, by introducing fractional derivatives
		in the constitutive equations of a one-dimensional compressible fluid,
		we embrace a vast phenomenology, including subdiffusive and superdiffusive processes.
		Moreover, with our mathematical hypotesis, processes like classical diffusions,
		can be considered as special cases of a more
		general and complex phenomenology that can be treated with fractional calculus.
		This means that in some physical problems, a diffusive 
		equation corresponds to a microscopic fractional process of migrating fluid particles with memory.
		From a mathematical point of view, we study different cases: when the semigroup property is valid
		for the Caputo fractional derivative, we derive a fractional diffusion equation by decoupling the
		related system of equations;
		otherwise we obtain partial differential equations with fractional sequential derivatives.
		The organisation of this paper is the following: in Section \ref{doc}
		we recall some basic concepts and properties related to fractional calculus;
		in Section \ref{doc2}, our model of propagation of a small perturbation 
		with time fractional derivative is discussed; in Section \ref{doc3}
		we provide a mathematical discussion on the range of variation of the characteristic indices
		of differentiation; Section \ref{stapler} concerns the solution of initial
		value problems involving fractional	sequential derivatives. This is particularly important as
		they arise naturally for some values of the characteristic indices of fractionality;
		finally, in Section \ref{doc5} we give a complete
		discussion of the results in relation to the previous research on generalised sound propagation models.

	\section{Basic concepts and theorems}

		\label{doc}
	    In this section we recall some definitions and basic results related to fractional calculus.

	 	\begin{Definition}
			Let $\alpha \in \mathbb{R}^{+}$. The Riemann--Liouville fractional integral is defined as
			\begin{equation}
				J^{\alpha}_t f(t) = \frac{1}{\Gamma(\alpha)}\int_0^{t}(t-\tau)^{\alpha-1}f(\tau) \mathrm \, d\tau,
				\qquad t>0.
			\end{equation}
		\end{Definition}
		Note that this operator is such that $J^0_t f(t)= f(t)$.
		Moreover, it satisfies the semigroup property, i.e. $J^{\alpha}J^{\beta} f(t) = J^{\alpha+\beta}f(t)$,
		$\alpha>0$, $\beta>0$.
		
		Different definitions of fractional derivatives have been introduced
		and studied in the literature (see for example \cite{Podl} ot \cite{kilbas}). For the sake
		of our analysis and because it proved to be of practical application, we make use of
		the Caputo fractional derivative.

		\begin{Definition}
			\label{Caputo}
			Let $m-1 < \alpha\leq m$, with $m\in \mathbb{N}$, the Caputo fractional derivative is defined by
			\begin{equation}
				D_t^{\alpha}f(t)= J^{m-\alpha}_t D_t^m f(t)= 
				\begin{cases}
					\frac{1}{\Gamma(m-\alpha)}\int_0^{t}(t-
					\tau)^{m-\alpha-1}f^{(m)} (\tau) \, \mathrm d\tau, & \alpha \ne m, \\
					f^{(m)}(t), & \alpha = m,
				\end{cases}
			\end{equation}
			where $f^{(m)}(t)=\frac{d^m f(t)}{dt^m}$ is the ordinary derivative of integer order m.	
		\end{Definition} 
		It is clear by Definition \ref{Caputo} that the fractional derivative is a pseudo-differential operator
		given by the convolution
		of the ordinary derivative with a power law kernel. Therefore the reason why fractional derivatives introduce
		a memory formalism becomes evident. 
		It is simple to prove the following properties of fractional derivatives and integrals (see e.g.\ \cite{Podl}):
		\begin{align}
			&D^{\alpha}_t J^{\alpha}_t f(t)= f(t), \qquad \alpha> 0,\\
			&J^{\alpha}_t D^{\alpha}_t f(t)= f(t)-\sum_{k=0}^{m-1}f^{(k)}(0)\frac{t^k}{k!}, \qquad \alpha>0, \: t>0,\\
			&D^{\alpha}_t t^{\beta}= \frac{\Gamma(\beta+1)}{\Gamma(\beta-\alpha+1)}t^{\beta-\alpha}
			\qquad \alpha>0, \: \beta\in\mathbb{R}, \: t>0.
		\end{align}

		Finally, a key point for the following discussion is the limit of validity of the semigroup property 
		for the Caputo fractional derivative.
		In its recent book, Diethelm (see \cite{Diethelm}, pag.56) gave a sufficient condition for the validity 
		of the semigroup property for the Caputo fractional derivative. Let us recall it.

		\begin{Theorem}[Law of exponents]
			\label{semi}
			Let $f(t) \in C^k[0,x]$, for some $x>0$ and some $k\in \mathbb{N}$. Moreover, let $\alpha, \beta >0$ 
			be such that there exists some $l\in \mathbb{N}$
			with $l\leq k$ and $\alpha, \alpha+\beta\in[l-1, l]$. Then
			\begin{equation}
				D^{\beta}_t D^{\alpha}_t f(t) = D^{\alpha+\beta}_t f(t).
			\end{equation} 
		\end{Theorem}
		The proof of this theorem is very simple. It proceeds by considering the definition of the Caputo derivative in
		relation to the Riemann--Liouville fractional integral.
		This theorem highlights a constraint on the applicability of the semigroup both with respect to the request
		of smoothness of the function and with respect to the ranges
		of the real orders of differentiation $\alpha$ and $\beta$. This means, for example, that,
		if $\alpha \in (0,1]$, then the law of exponents
		is applicable if $\beta \in [0, 1-\alpha)$ and $f(t)\in C^k$, with $k \geq 1$. 
		Here we have also to notice that in most cases the law of exponents is not applicable for
		fractional Caputo derivatives, but anyhow there are different techniques
		to handle sequential fractional derivatives (see for example \cite{Podl}). 
		It is clear that in cases where the semigroup is still valid for the Caputo derivative, the 
		fractional equations are more manageable and, what is most 
		important, the construction of a Cauchy problem demands conditions with meaningful physical sense. In the 
		following discussion therefore, we develop our analysis
		in the cases in which the semigroup property is valid, neglecting the analysis of the fractional sequential 
		case for which the requested initial conditions have not 
		a clear physical interpretation.
	
	\section{Propagation of a small perturbation by means of fractional calculus modelling}

		\label{doc2}
		Let us consider a one-dimensional compressible fluid, initially at rest, with an homogeneous 
		constant density $\rho_0$. We study the effect of a small
		perturbation in the density $\rho'(z,t)$. We thus deal with a perturbative theory with
		\begin{equation}
			\rho(z,t) = \rho_0 + \rho'(z,t), \qquad w(z,t) = w'(z,t), \qquad z \in \mathbb{R}, \: t\geq 0,
		\end{equation}
		where $w(z,t)$ is the velocity field and where $\rho'(z,t)$ and $w'(z,t)$ are first order terms.
		This is a simple way to obtain the sound wave propagation equation in a compressible fluid. 
		We will develop this classical model replacing the ordinary derivative with the fractional Caputo derivative.
		
		The equations of continuity and velocity of the fluid, linearised at the first order, are

		\begin{align}
			\label{dist}
			\begin{cases}
				\partial^{\alpha}_t \rho'(z,t)+ \rho_0 \partial_z w'(z,t) = 0, & \alpha>0, \: \beta>0,\: t\ge 0, \\
				\rho_0 \partial^{\beta}_t w'(z,t) + c_s^2\partial_z \rho'(z,t)= 0,
			\end{cases}
	  	\end{align}
		where $c_s^2=\partial p'/\partial \rho'$ is the sound velocity in the fluid. 
		This is a coupled linear system of fractional equations.
		Note that we have used the symbol $\partial^{\alpha}_t$ for the Caputo 
		fractional partial derivative and in general the real order $\alpha$ and $\beta$
		can be different. In the following we discuss the range of variability of the fractional orders in $\mathbb{R}$.
		
		Now, applying the fractional derivative of order $\beta$
		to the first equation, and differentiating the second equation with respect to $z$, we obtain

		\begin{align}
			\begin{cases}
				\partial^{\beta}_t\partial^{\alpha}_t \rho'(z,t)+ \rho_0 \partial_z\partial^{\beta}_t w'(z,t) = 0, \\ 
				\partial_z \partial^{\beta}_t w'(z,t) +\frac{c_s^2}{\rho_0}\partial_{zz} \rho'(z,t) =0,
			\end{cases}
	  	\end{align}
		and finally
		\begin{equation}
			\label{seq}
			\partial^{\beta}_t \partial^{\alpha}_t \rho'(z,t)- c_s^2 \partial_{zz}\rho'(z,t)=0,
			\qquad \alpha>0, \: \beta>0, \: t\ge 0.
		\end{equation}
		In general this is a partial differential equation with a time fractional sequential derivative (see \cite{Podl}).
		We are interested in the cases in which the law of exponents (Theorem \ref{semi})
		is satisfied and the equation \eqref{seq} therefore becomes	
		\begin{equation}
			\partial^{\alpha+\beta}_t \rho'(z,t)- c_s^2 \partial_{zz}\rho'(z,t)=0, \qquad t \ge 0, \: z \in \mathbb{R},
		\end{equation}
		i.e.\ a fractional diffusion-wave equation (when $\alpha + \beta \in (0,2]$),
		extensively studied in literature (see for example \cite{Main} 
		for the fundamental solution). In fact we wish to arrive at a meaningful equation that can be treated 
		considering initial conditions with physical meaning,
		differently from the case usually treated with fractional sequential derivatives. 
		In the next section we thus discuss the cases in which this coupled system leads to a single 
		fractional diffusion-wave equation.
		Note also that with a symmetric manipulation we find the following equation for the velocity field
		$w'(z,t)$:
		\begin{equation} 
			\partial^{\alpha}_t \partial^{\beta}_t w'(z,t)- c_s^2 \partial_{zz}w'(z,t)=0, \qquad t \ge 0, \: z \in
			\mathbb{R}.
		\end{equation}
		 
	\section{Coupled systems of fractional equations related to diffusion and subdiffusion}

		\label{doc3}
	 	In the previous section, we showed that our coupled system of fractional equations on 
	 	velocity and density fields of a one-dimensional
		compressible fluid, is related to the fractional diffusion-wave equation, under 
		the hypothesis that the semigroup property holds.
		Here we discuss in detail the cases in which we can apply the law of exponents. 
		This section is devoted to mathematical analysis; we will return
		to physical meaning in the last section.
		
\begin{figure}
			\centering
			\includegraphics[scale=.35]{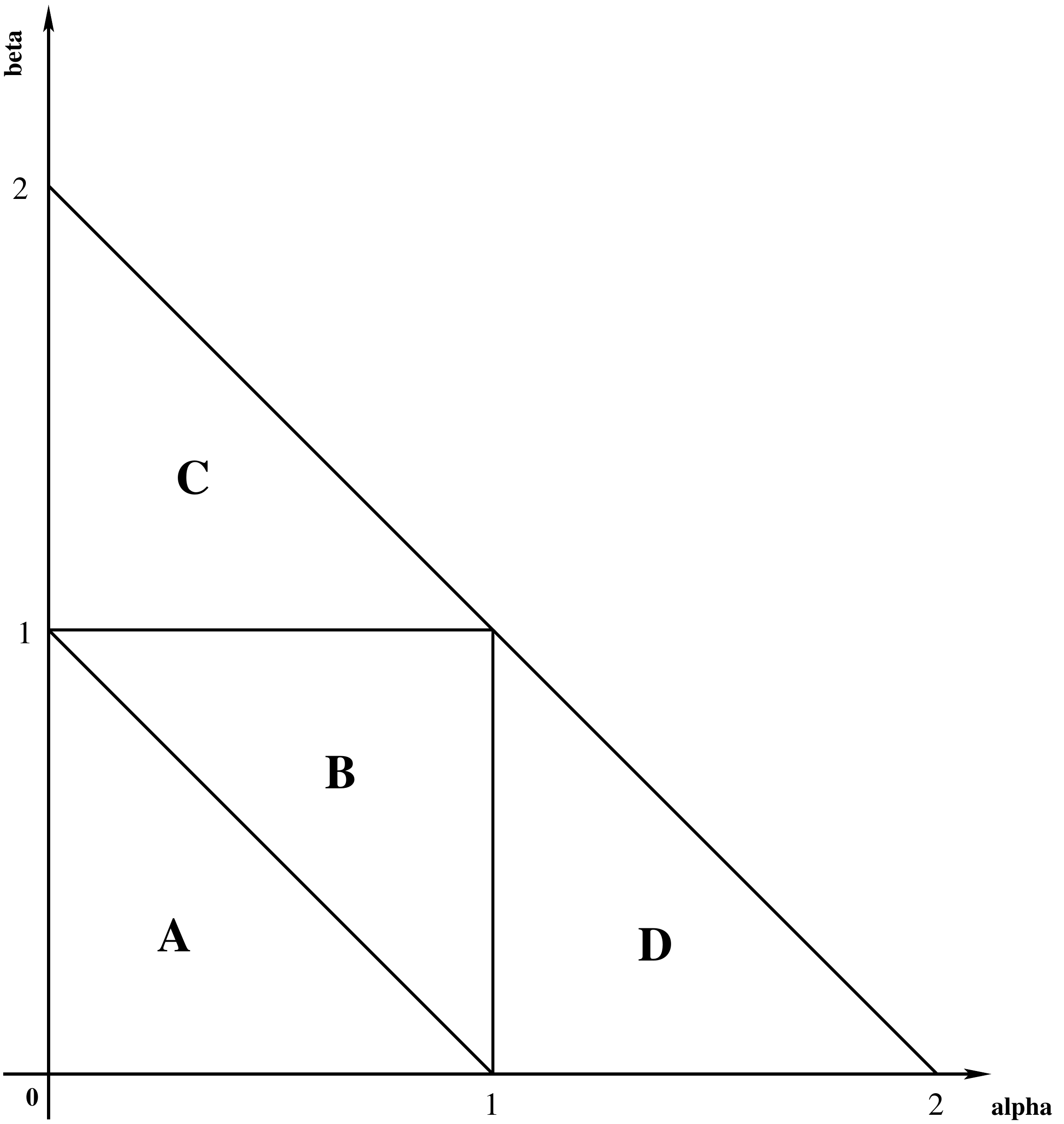}
			\caption{\label{tione3}The four characteristic regions in the $\alpha\beta$-plane.}
		\end{figure}

		We first consider equation \eqref{seq} with the constraints
		\begin{equation}
			\label{ciao}
			\alpha \in (0,1), \qquad \beta \in (0, 1-\alpha],
		\end{equation}
		i.e. the pair of parameters $(\alpha,\beta)$ varies in region $A$ shown in Fig.\ \ref{tione3}.
		Note that the constraints \eqref{ciao} are equivalent to
		\begin{align}
			\label{ciao2}
			\beta \in (0,1), \quad \alpha \in (0, 1-\beta].
		\end{align}
		
		When \eqref{ciao} or \eqref{ciao2} hold, for the Caputo fractional derivative a semigroup
		rule is valid (see Theorem \ref{semi}), then we can rewrite \eqref{seq} as
		\begin{align}
			\label{statua}
			\partial^{\alpha + \beta}_t \rho'(z,t) - c_s^2 \partial_{zz} \rho'(z,t) = 0
			\qquad \overset{\alpha+\beta=\gamma}{\Leftrightarrow}
			\qquad \partial^\gamma_t \rho'(z,t) - c_s^2 \partial_{zz}\rho'(z,t) = 0.
		\end{align}
		Equation \eqref{statua} is clearly a time-fractional diffusion equation (see \cite{Podl}, page 296).
		Conditions \eqref{ciao} implies that $\gamma \in (0,1]$, therefore the regime is subdiffusive
		or at most diffusive. The diffusive regime (corresponding to the classical heat equation) is reached
		for each pair $(\alpha,\beta)$ for which $\beta = 1-\alpha$ (see Fig.\ \ref{tione3}).
		
		We are now going to study the possibility of arriving at a fractional diffusion equation 
		displaying a superdiffusive behaviour. 
		In this case, in order to maintain at least partially the semigroup property,
		we consider the following constraints for the fractional indices $\alpha$ and $\beta$:
		\begin{align}
			\label{ciao3}
			\alpha \in (1,2), \qquad \beta \in (0, 2-\alpha].
		\end{align}
		In Fig.\ \ref{tione3} the region $D$ of the $\alpha \beta$-plane in which \eqref{ciao3} hold is shown.
		
		Equation \eqref{statua} is easily retrieved. Clearly, in this case $\gamma = \alpha+\beta \in (1,2]$
		thus furnishing the system with a superdiffusive regime. Furthermore the related coupled fractional
		system is \eqref{dist} where \eqref{ciao3} hold.

		\begin{Remark}
			\label{samsung}
			If conditions \eqref{ciao3} hold, the equation related to the velocity field arising from \eqref{dist} 
			can be determined as
			\begin{align}
				\label{statua2}
				\partial^\alpha_t \partial^\beta_t w'(z,t) = c_s^2 \partial_{zz} w'(z,t), \qquad t \ge 0,
				\: z \in \mathbb{R}.
			\end{align}
			Here the semigroup property is no longer valid.
		\end{Remark}
		
		\begin{Remark}
			Considerations similar to those in Remark \ref{samsung} in the case of
			\begin{align}
				\label{ciao4}
				\beta \in (1,2), \qquad \alpha \in (0, 2-\beta],
			\end{align}
			lead to the following results. The semigroup validity region is pictured in Fig.\ \ref{tione3}, region
			$C$.
			In this case we obtain
			\begin{align}
				\label{system}
				\partial^\gamma_t w'(z,t) = c^2_s \partial_{zz} w'(z,t), \qquad t \ge 0, \: z \in \mathbb{R},
			\end{align}
			and
			\begin{align}
				\label{system2}
				\partial^\beta_t \partial^\alpha_t \rho'(z,t) = c^2_s \partial_{zz} \rho'(z,t)
				\qquad t \ge 0, \: z \in \mathbb{R}.
			\end{align}
		\end{Remark}
		
		Concluding we note that the cases for which the semigroup holds, are not 
		exactly a generalisation of the classical case ($\alpha=1$, $\beta=1$).
		However we will give in the last section a physical 
		interpretation of these equations in relation to the original
		framework, i.e.\ equations of velocity and mass conservation of one-dimensional compressible fluid. 
		
		\begin{Remark}
			Here we remark that solutions to both uncoupled equations
			\eqref{statua} and \eqref{statua2}
			(or the related \eqref{system} and \eqref{system2})
			can be determined by first solving the subdiffusion type equation and then by arriving
			at a solution of the sequential derivative type equation by means of the relations of the coupled
			system \eqref{dist}.
			Indeed we have for example that the solution to the Cauchy problem
			\begin{align}
				\begin{cases}
					\partial^\gamma_t \rho'(z,t) = c^2_s \partial_{zz} \rho'(z,t), & t\ge 0, \: z \in \mathbb{R}, \\
					\rho'(z,0) = \delta(z), \\
					\partial_t \rho'(z,t)|_{t=0} = 0,
				\end{cases}
			\end{align}
			is well-known and reads
			\begin{align}
				\rho'(z,t) = \frac{1}{2c^2_s t^{\gamma/2}} W_{-\gamma/2, 1-\gamma/2}
				\left( -\frac{|z|}{c^2_s t^{\gamma/2}} \right), \qquad t\ge 0, \: z \in \mathbb{R},
			\end{align}
			where
			\begin{align}
				W_{\kappa,\eta}(y) = \sum_{r=0}^\infty \frac{y^r}{r!\Gamma(\kappa r+\eta)},
				\qquad \kappa>-1, \: y \in \mathbb{C}
			\end{align}
			is the Wright function \citep[page 54]{kilbas}. In turn, with the initial condition
			$w'(z,0) = f(z)$, we can write
			\begin{align}
				& \partial^\beta_t w'(z,t) = - \frac{c^2_s}{\rho_0} \partial_z \rho'(z,t) \\
				& \Leftrightarrow \quad J^\beta_t \partial^\beta_t
				w'(z,t) = - \frac{c^2_s}{\rho_0} J^\beta_t \partial_z \rho'(z,t) \notag \\
				& \Leftrightarrow \quad w'(z,t) = f(z) - \frac{c^2_s}{\rho_0} J^\beta_t \partial_z
				\rho'(z,t) \notag.
			\end{align}
		\end{Remark}
		
	\section{Coupling for PDEs involving sequential fractional derivatives}
	
		\label{stapler}
		In this section we present a possible way of addressing the problem of solving Cauchy problems
		involving sequential fractional derivatives. We consider the region of the $\alpha \beta$-plane
		defined by the constraint $(\alpha,\beta) \in B$ (see Fig.\ \ref{tione3}). We seek to solve the following
		initial value problem:
		\begin{align}
			\label{inter}
			\begin{cases}
				\partial_t^\alpha \partial_t^\beta f(z,t) = \lambda \partial_{zz}f(z,t), & t \ge 0, \:
				z \in \mathbb{R}, \\
				f(z,0) = g(z), \\
				\partial_t^\beta f(z,t)|_{t=0} = \bar{g} (z).
			\end{cases}
		\end{align}
		We underline that we have requested specific initial conditions involving fractional derivatives.
		This choice will be clear in the following discussion. Here we recall that this kind of initial conditions has not 
		a broadly accepted physical meaning,
		but some suggestions on its interpretation are already present in the literature
		(see \cite{Hey}). In this section however, we are more interested in outlining a possible
		way to solve partial differential equations with fractional sequential equations more than
		in defining a precise physical meaning.
		 
		In order to solve \eqref{inter}, we consider an auxiliary function $\varphi(z,t)$, solution to
		\begin{align}
			\label{inter2}
			\begin{cases}
				\partial_t^\beta \partial_t^\alpha \varphi(z,t) = \lambda \partial_{zz}\varphi(z,t), & t \ge 0, \:
				z \in \mathbb{R}, \\
				\varphi(z,0) = h(z), \\
				\partial_t^\alpha \varphi(z,t)|_{t=0} = \bar{h} (z).
			\end{cases}
		\end{align}
		By means of \eqref{inter} and \eqref{inter2}, we can construct the coupled system
		\begin{align}
			\label{44}
			\begin{cases}
				\partial_t^\beta f(z,t) + \frac{\lambda}{\kappa} \partial_z \varphi(z,t) = 0, \\
				\partial_t^\alpha \varphi(z,t) + \kappa \partial_z f(z,t) = 0,
			\end{cases}
		\end{align}
		with the initial conditions $f(z,0)=g(z)$ and $\varphi(z,0)=h(z)$. Notice that functions
		$g(z)$ and $h(z)$ are linked together. Indeed from \eqref{44} we have that
		\begin{align}
			\label{lab}
			& \partial_t^\beta f(z,t) |_{t=0} = -\frac{\lambda}{\kappa} \partial_z \varphi(z,0) \\
			& \Leftrightarrow \quad \bar{g}(z) = -\frac{\lambda}{\kappa} \partial_z \varphi(z,0). \notag
		\end{align}
		By inserting the initial condition of \eqref{inter2} into \eqref{lab} we immediately obtain
		\begin{align}
			\bar{g}(z) = -\frac{\lambda}{\kappa} D_z h(z).
		\end{align}
		Similarly, from \eqref{44},
		\begin{align}
			\partial_t^\alpha \varphi(z,t)|_{t=0} = -\kappa \partial_z f(z,0).
		\end{align}
		Considering \eqref{inter2} and \eqref{inter} we arrive at
		\begin{align}
			\bar{h}(z) & = -\kappa \partial_z f(z,0)
			= -\kappa D_z g(z).
		\end{align}
		
		Summarising, the coupled Cauchy problems becomes
		\begin{align}
			\label{allu}
			\begin{cases}
				\partial_t^\beta f(z,t) + \frac{\lambda}{\kappa} \partial_z \varphi(z,t) = 0, \\
				\partial_t^\alpha \varphi(z,t) + \kappa \partial_z f(z,t) = 0, \\
				f(z,0) = g(z), \qquad \bar{g}(z) = -\frac{\lambda}{\kappa} D_z h(z), \\
				\varphi(z,0) = h(z), \qquad \bar{h}(z) = -\kappa D_z g(z),
			\end{cases}
		\end{align}
		$t \ge 0$, $z \in \mathbb{R}$. Now, the application of the Fourier transform
		\begin{align}
			\hat{\zeta}(\omega,t) = \int_{-\infty}^\infty e^{i\omega x} \zeta(x,t) \, \mathrm dx
		\end{align}
		to \eqref{allu} leads to
		\begin{align}
			\label{alluf}
			\begin{cases}
				\partial_t^\beta \hat{f}(\omega,t) = - \frac{\lambda}{\kappa} i\omega \hat{\varphi}(\omega,t), \\
				\partial_t^\alpha \hat{\varphi}(\omega,t) = - \kappa i \omega \hat{f}(\omega,t), \\
				\hat{f}(\omega,0) = \hat{g}(\omega), \qquad \hat{\bar{g}}(\omega)
				= -\frac{\lambda}{\kappa} i \omega \hat{h}(\omega), \\
				\hat{\varphi}(\omega,0) = \hat{h}(\omega), \qquad \hat{\bar{h}}(\omega)
				= -\kappa i \omega \hat{g}(\omega).
			\end{cases}
		\end{align}
		In turn, by passing to the Laplace transform
		\begin{align}
			\tilde{\zeta}(x,s) = \int_0^\infty e^{-st} \zeta(x,t) \, \mathrm dt
		\end{align}
		we have
		\begin{align}
			\begin{cases}
				s^\beta \tilde{\hat{f}}(\omega,s) -s^{\beta-1} \hat{g}(\omega)
				= -\frac{\lambda}{\kappa} i \omega \tilde{\hat{\varphi}}(\omega,s), \\
				s^\alpha \tilde{\hat{\varphi}}(\omega,s) - s^{\alpha-1} \hat{h}(\omega)
				= -\kappa i \omega \tilde{\hat{f}}(\omega,s).
			\end{cases}
		\end{align}
		By substitution we can write that
		\begin{align}
			& s^\beta \tilde{\hat{f}}(\omega,s) - s^{\beta-1} \hat{g}(\omega) =
			-\frac{\lambda}{\kappa} i \omega \left( s^{\alpha-1} \hat{h}(\omega) -
			\kappa i \omega \tilde{\hat{f}}(\omega,s) \right)\frac{1}{s^\alpha} \\
			& \Leftrightarrow \quad s^{\alpha+\beta} \tilde{\hat{f}} (\omega,s) - s^{\alpha+\beta-1}
			\hat{g}(\omega) = -\frac{\lambda}{\kappa} i \omega s^{\alpha-1}
			\hat{h}(\omega) - \lambda \omega^2 \tilde{\hat{f}}(\omega,s) \notag \\
			& \hspace{-.35cm} \overset{\text{from} \eqref{alluf}}{\Leftrightarrow} \quad
			s^{\alpha+\beta} \tilde{\hat{f}} (\omega,s) - s^{\alpha+\beta-1}
			\hat{g}(\omega) = s^{\alpha-1} \hat{\bar{g}}(\omega) - \lambda \omega^2 \tilde{\hat{f}}(\omega,s) \notag \\
			& \Leftrightarrow \quad \tilde{\hat{f}}(\omega,s) =
			\hat{g}(\omega) \frac{s^{\alpha+\beta-1}}{s^{\alpha+\beta} + \lambda \omega^2}
			+ \hat{\bar{g}}(\omega) \frac{s^{\alpha-1}}{s^{\alpha+\beta}+\lambda \omega^2}. \notag
		\end{align}
		We now make use of the following Laplace transform:
		(see \citet{haubold}, formula (2.2.26))
		\begin{align}
			\int_0^\infty e^{-st} t^{\gamma -1} E_{\eta,\gamma} (a t^\eta) \, \mathrm dt = \frac{s^{\eta-\gamma}}{s^\eta - a},
		\end{align}
		where
		\begin{align}
			E_{\eta,\gamma}(y) = \sum_{r=0}^\infty \frac{y^r}{\Gamma(\eta r+ \gamma)}, \qquad y \in \mathbb{R},
		\end{align}
		is the Mittag--Leffler function. Therefore, by inverting the Laplace transform we obtain
		\begin{align}
			\hat{f}(\omega,t) = \hat{g}(\omega) E_{\alpha+\beta,1}(-\lambda \omega^2 t^{\alpha+\beta})
			+ \hat{\bar{g}}(\omega) t^\beta E_{\alpha+\beta,\beta+1} (-\lambda \omega^2 t^{\alpha+\beta}).
		\end{align}
		The inversion of the Fourier transform is now only a matter of an application of relations
		(14.7) and (14.8), page 27 of \citet{hau} which leads to an explicit formula involving Fox functions.
		\begin{align}
			f(z,t) = {} & g(z) \, \ast \, \left( \frac{1}{|z|} H^{1,0}_{1,1} \left[ \left.
			\frac{|z|^2}{\lambda t^{\alpha+\beta}} \right|
			\begin{array}{l}
				(1,\alpha+\beta) \\
				(1,2)
			\end{array}
			\right] \right) \\
			& + t^\beta \hat{g}(z) \, \ast \, \left( \frac{1}{|z|} H^{1,0}_{1,1} \left[ \left.
			\frac{|z|^2}{\lambda t^{\alpha+\beta}} \right|
			\begin{array}{l}
				(1+\beta,\alpha+\beta) \\
				(1,2)
			\end{array}
			\right] \right), \qquad t \ge 0, \: z \in \mathbb{R}, \notag
		\end{align}
		where the symbol $\ast$ represents the Fourier convolution.
		From the definition of the Fox function (see e.g.\ \citet{kilbas}) we can write
		\begin{align}
			f(z,t) = {} & g(z) \, \ast \, \left( \frac{1}{|z|} \frac{1}{2 \pi i} \int_C \frac{\Gamma(1+2\vartheta)}{
			\Gamma(1+(\alpha+\beta)\vartheta)} \left( \frac{|z|^2}{\lambda t^{\alpha+\beta}}
			\right)^{-\vartheta} \mathrm d \vartheta \right) \\
			& + t^\beta \bar{g}(z) \, \ast \, \left( \frac{1}{|z|} \frac{1}{2 \pi i} \int_C \frac{\Gamma(1+2\vartheta)}{
			\Gamma(1+\beta +(\alpha+\beta)\vartheta)} \left( \frac{|z|^2}{\lambda t^{\alpha+\beta}}
			\right)^{-\vartheta} \mathrm d \vartheta \right) \notag \\
			= {} & g(z) \, \ast \, \left( \frac{1}{\sqrt{\lambda t^{\alpha+\beta}}} \frac{1}{2\pi i}
			\int_C \frac{\Gamma(1+2\vartheta)}{\Gamma(1+(\alpha+\beta)\vartheta)}
			\left( \frac{|z|}{\sqrt{\lambda t^{\alpha+\beta}}} \right)^{-2\vartheta-1} \mathrm d \vartheta \right) \notag \\
			& + t^\beta \bar{g}(z) \, \ast \, \left( \frac{1}{\sqrt{\lambda t^{\alpha+\beta}}} \frac{1}{2\pi i}
			\int_C \frac{\Gamma(1+2\vartheta)}{\Gamma(1+\beta+(\alpha+\beta)\vartheta)}
			\left( \frac{|z|}{\sqrt{\lambda t^{\alpha+\beta}}} \right)^{-2\vartheta-1} \mathrm d \vartheta \right) \notag \\
			= {} & g(z) \, \ast \, \left( \frac{1}{2\sqrt{\lambda t^{\alpha+\beta}}} \frac{1}{2\pi i}
			\int_C \frac{\Gamma(\tau)}{\Gamma\left( 1+(\alpha+\beta)\frac{\tau-1}{2} \right)}
			\left( \frac{|z|}{\sqrt{\lambda t^{\alpha+\beta}}} \right)^{-\tau} \mathrm d \tau \right) \notag \\
			& + t^\beta \bar{g}(z) \, \ast \, \left( \frac{1}{2\sqrt{\lambda t^{\alpha+\beta}}} \frac{1}{2\pi i}
			\int_C \frac{\Gamma(\tau)}{\Gamma\left( 1+\beta+(\alpha+\beta)\frac{\tau-1}{2} \right)}
			\left( \frac{|z|}{\sqrt{\lambda t^{\alpha+\beta}}} \right)^{-\tau} \mathrm d \tau \right) \notag \\
			= {} & g(z) \, \ast \, \left( \frac{1}{2\sqrt{\lambda t^{\alpha+\beta}}} \frac{1}{2\pi i}
			\int_C \frac{\Gamma(\tau)}{\Gamma\left( 1-\frac{\alpha+\beta}{2} + \frac{\alpha+\beta}{2} \tau \right)}
			\left( \frac{|z|}{\sqrt{\lambda t^{\alpha+\beta}}} \right)^{-\tau} \mathrm d \tau \right) \notag \\
			& + t^\beta \bar{g}(z) \, \ast \, \left( \frac{1}{2\sqrt{\lambda t^{\alpha+\beta}}} \frac{1}{2\pi i}
			\int_C \frac{\Gamma(\tau)}{\Gamma\left( 1+\beta - \frac{\alpha+\beta}{2}
			+\frac{\alpha+\beta}{2} \tau \right)}
			\left( \frac{|z|}{\sqrt{\lambda t^{\alpha+\beta}}} \right)^{-\tau} \mathrm d \tau \right). \notag
		\end{align}
		Furthermore, by means of the Mellin--Barnes representation of the Wright function \citep[formula (1.11.3)]{kilbas}
		we finally obtain
		\begin{align}
			\label{treno}
			f(z,t) = {} & g(z) \, \ast \, \left( \frac{1}{2\sqrt{\lambda t^{\alpha+\beta}}} W_{-\frac{\alpha+\beta}{2},
			1-\frac{\alpha+\beta}{2}} \left( -\frac{|z|}{\sqrt{\lambda t^{\alpha+\beta}}} \right) \right) \\
			& + t^\beta \bar{g}(z) \, \ast \, \left( \frac{1}{2\sqrt{\lambda t^{\alpha+\beta}}} W_{-\frac{\alpha+\beta}{2},
			1+\beta-\frac{\alpha+\beta}{2}} \left( -\frac{|z|}{\sqrt{\lambda t^{\alpha+\beta}}} \right) \right),
			\qquad t \ge 0, \: z \in \mathbb{R}. \notag
		\end{align}
		
		\begin{Remark}
			Note that \eqref{treno} when $\bar{g}(z) = 0$ concides with the solution to the
			superdiffusive fractional diffusion equation. The reader can consult \citet{kilbas}, Section 6.1.2.
			In this case \eqref{treno} shows also that the semigroup property is satisfied.
		\end{Remark}
		
		\begin{Remark}
			We note that in the case $(\alpha,\beta) \in D$ (or correspondingly, $(\alpha,\beta) \in C$), i.e.\
			when the problems
			\begin{align}
				\label{interb}
				\begin{cases}
					\partial_t^\alpha \partial_t^\beta f(z,t) = \lambda \partial_{zz}f(z,t), & t \ge 0, \:
					z \in \mathbb{R}, \\
					f(z,0) = g(z), \\
					\partial_t^\beta f(z,t)|_{t=0} = \bar{g} (z),
				\end{cases}
			\end{align}
			and
			\begin{align}
				\label{inter2b}
				\begin{cases}
					\partial_t^{\alpha + \beta} \varphi(z,t) = \lambda \partial_{zz}\varphi(z,t), & t \ge 0, \:
					z \in \mathbb{R}, \\
					\varphi(z,0) = h(z),
				\end{cases}
			\end{align}
			are taken into consideration, a reasoning similar to that of this section, leads to the same
			solution \eqref{treno} for $f(z,t)$ and to the solution to the well-known fractional diffusion-wave
			equation for function $\varphi(z,t)$, as it should be.
		\end{Remark}

	\section{Discussion}
	
		\label{doc5}
		Here we discuss the main result of this paper: finding diffusive processes from the 
		fractional generalisation of density and velocity equations
		of a simple one-dimensional compressible fluid.
		
		From a physical point of view we tried to understand the consequence of using a wide generalisation of
		the constitutive equations 
		on sound wave propagation in order to take into account memory effects. In this way we want to show that 
		fractional calculus proves to be a useful analytic 
		tool to describe a general class of phenomena beyond classical memoryless processes. 
		Moreover, our aim is to build a complete theoretical framework to justify the empirical modification of the wave
		equation discussed in literature. Regarding this, for example, we recall the recent paper of \citet{Prieur} that
		relates non-integer power law absorption of sound waves to fractional equations; they justify the modification of
		the wave equation starting from the constitutive equations rather than considering empirical modifications due to
		observed attenuation and dispersion power laws.		
		We also suggest an interpretation within the 
		framework of fluid mechanics, recalling
		some relevant recent results on fractional modelling in this field.
		In their work Fellah et al.\ \cite{Fellah2} suggested the following constitutive equations of
		acoustic waves propagations in porous media:

		\begin{align}
			\label{Fel}
			\begin{cases}
				\rho\tilde{\alpha}(t)\ast \partial_t w = -\phi \partial_x p, \\
				\frac{\phi\tilde{\beta}(t)}{K_a}\ast \partial_t p = - \partial_x w,
			\end{cases}
		\end{align}
		where $K_a$ is the bulkmodulus of air, $\phi$ is the porosity, $\rho$ the density, 
		$p$ is the acoustic pressure, $\tilde{\alpha}(t)$ and 
		$\tilde{\beta}(t)$ are the dynamic tortuosity and compressibility of the air; clearly the 
		symbol $\ast$ denotes the time convolution
		\begin{equation}
			(f\ast g)(t)=\int_{0}^{t}f(t-t')g(t')dt'.
		\end{equation}
		We refer to the original paper for an analysis of the physical model. What is most 
		relevant is that they used fractional calculus
		to take into account time variability of $\tilde{\alpha}$ and $\tilde{\beta}$. The 
		convolution with dynamical coefficients clearly implies 
		a physical global delay in the signal propagation and the presence of memory. 
		Also the experimental variability of dynamical coefficients
		suggests to use the time-fractional derivative in place of the convolution considered in \eqref{Fel}. 
		This fractional model for acoustic wave propagation was also experimentally validated in \cite{Fellah3}. 
		In this framework our model can be simply adapted. In fact, following \cite{Fellah2}, \eqref{Fel} can
		be written as a coupled system of fractional equations
		on pressure and velocity fields in a one-dimensional model. 
		The indices related to the order of differentiation are free parameters to be fit by means of experimental 
		findings on the variability of the characteristic coefficients in the medium of propagation.
		Furthermore there is no need to assume $\alpha = \beta$ 
		in our constitutive equations. Fractional calculus
		in this case emerges as a natural formalism from the formulation used in \cite{Fellah2}.
		This allows us to have a family of flexible processes through which we can deal
		in a general way with
		many different physical situations depending on the properties of the media of propagation.	
		Finally, considering the results biefly reviewed here, we understand that fractional wave equations are relevant
		in cases of propagation in media characterised by power law frequency-dependent physical coefficients.		
		
		From a mathematical point of view, we find two interesting results. First of all,
		when the semigroup property holds for 
		the Caputo derivative, from the constitutive equations of the model we obtain
		diffusive or subdiffusive behaviours. 
		Hence we can infer anomalous diffusive processes from fractional models of fluid such as that in \eqref{Fel}.
		On the other hand, we find a way to solve partial differential equations with fractional
		sequential derivatives based on analysing the related
		coupled system of equations.


\end{document}